\documentclass[12pt]{amsart}
\usepackage{geometry}                
\usepackage{graphicx}
\usepackage{amssymb}
\usepackage{epstopdf}
\title{Metaphors in systolic geometry}
\author{Larry Guth}

\begin{document}
\maketitle

This essay is about Gromov's systolic inequality.  We will discuss why the inequality
is difficult, and we will discuss several approaches to proving the inequality
based on analogies with other parts of geometry.  The essay does not contain
proofs.  It is supposed to be accessible to a broad audience.

The story of the systolic inequality begins in the 1940's with Loewner's theorem.

\newtheorem*{syst0}{Loewner's systolic inequality} 

\begin{syst0} (1949) If $(T^2, g)$ is 
a 2-dimensional torus with a Riemannian metric, then there is a non-contractible
curve $\gamma \subset (T^2, g)$ whose length obeys the inequality

$$ length (\gamma) \le C Area(T^2, g)^{1/2},$$

where $C = 2^{1/2} 3^{-1/4}$.

\end{syst0}

To get a sense of Loewner's theorem, let's look at some pictures of 2-dimensional tori
in $\mathbb{R}^3$.  

\begin{figure}[htbp]
     \centering
    \includegraphics[width=12cm]{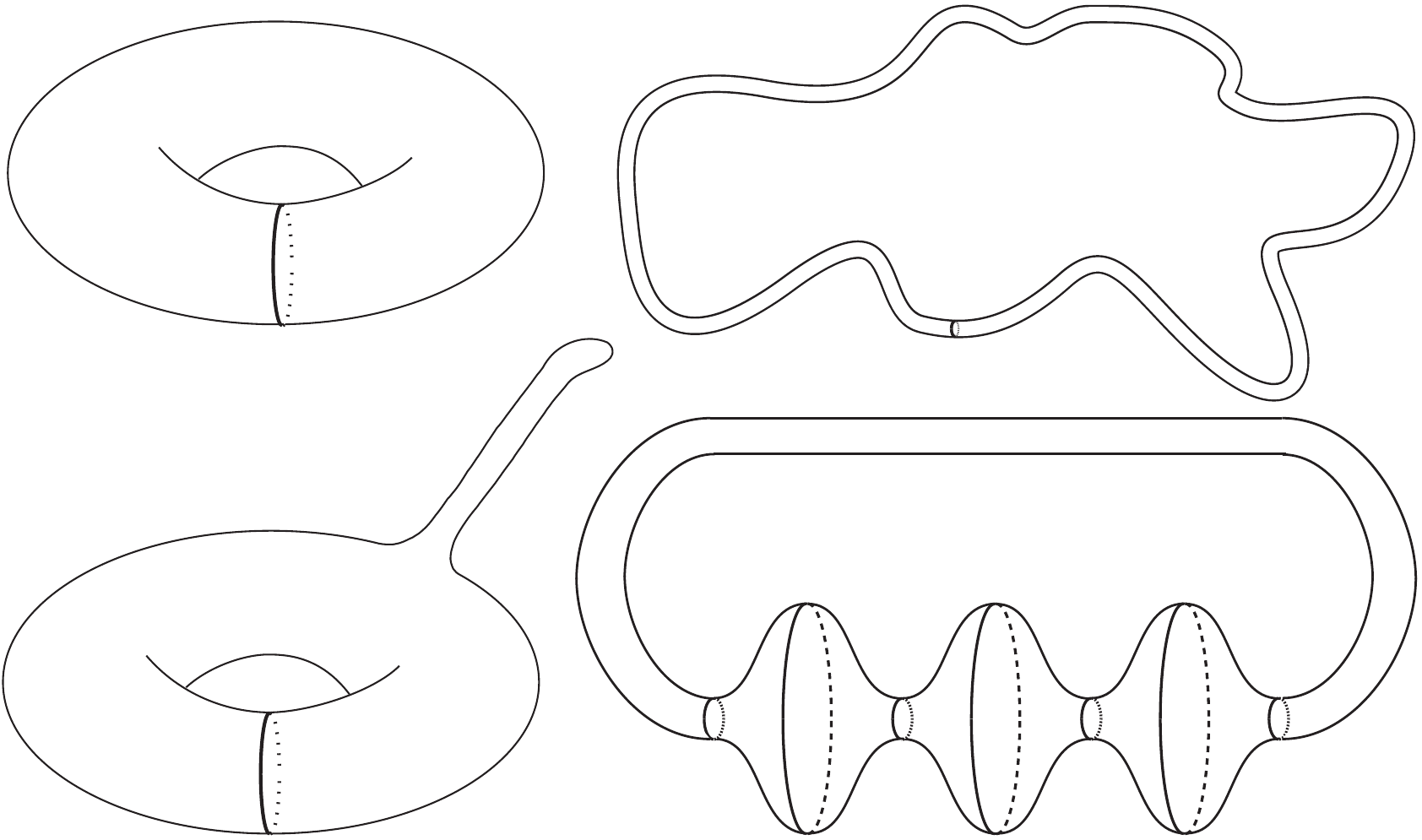}
    \caption{Pictures of tori}
    \label{cylinder}
\end{figure}

The curves shown in the pictures above are all non-contractible.  The length of
the shortest non-contractible curve on a Riemannian manifold is called its systole.

The first picture is supposed to show a torus of revolution, where we take the circle
of radius 1 around the point $(2,0)$ is the x-z plane and revolve it around the z-axis.
It has systole $2 \pi$ and area around 60,
and so it obeys the systolic inequality.  According to Loewner's theorem, there is nothing
we can do to dramatically increase the systole while keeping the area the same.
The second picture shows a long skinny torus.  When we make the torus skinnier and longer, the systole goes down and the area stays
about the same.  The third picture shows a torus with a long thin spike coming out of it.  
When we add a long thin spike to the torus, the systole doesn't change and the spike
adds to the area.  The fourth picture shows a ridged torus with some thick parts and some thin
parts.  When we put ridges in the surface of the torus, the systole only depends on the thinnest part and the thick
parts contribute heavily to the area.  These pictures are not a proof, but I think they make Loewner's
inequality sound plausible.

(Friendly challenge to the reader: can you think
of a torus with geometry radically different from the pictures above?)

Thirty years later, Gromov generalized Loewner's theorem to higher dimensions.

\newtheorem*{syst}{Gromov's systolic inequality for tori}

\begin{syst} (1983 \cite{G1}) If $(T^n, g)$ is an n-dimensional torus with a Riemannian metric, then the 
systole of $(T^n,g)$ is bounded in terms of its volume as follows.

$$ Sys(T^n,g) \le C_n Vol(T^n,g)^{1/n}. $$

\end{syst}

In his book {\it Metric Structures} \cite{G2}, Gromov reminisces about his work on the systolic inequality:
``Since the setting was so plain and transparent, I expected rather straightforward proofs
...  Having failed to find such a proof, I was inclined
to look for counterexamples, but..."  The statement of the theorem is extremely elementary
compared to other theorems in Riemannian geometry.  In spite of the plain and direct statement,
the theorem is difficult.  In particular, it's difficult to see how to approach the problem - how
to get started.

The systolic inequality for the 2-dimensional
torus was formulated and proven by Loewner, and a little later Besicovitch
gave a more elementary proof with a worse constant.  The two-dimensional proofs of
Loewner and Besicovitch do not generalize to three dimensions.  Gromov learned about
the problem in the late 60's from Burago, and it was popularized in the West by Berger.  Gromov
thought about it off and on during the 1970's and he devoted a chapter to it in the first
edition of \cite{G2}, published at the end of the 1970's.  At this point in time there was
still no good way to approach the systolic problem for the 3-dimensional torus.

In the early 80's, Gromov formulated several remarkable metaphors connecting the systolic
inequality to important ideas in other areas of geometry.  With the help
of these metaphors, he proved the systolic inequality.  We now have three independent 
proofs of the systolic inequality for the n-dimensional torus, each based on a different metaphor.
Each metaphor gives an approach to
proving the systolic inequality - a way to get started.

The goal of this essay is to explain Gromov's metaphors.  In doing that, I hope to describe the flavor of this branch of
geometry and put it into a broad context.  Gromov's metaphors connect the systolic
problem to the following areas:

1. General isoperimetric inequalities from geometric measure theory.  (Work of
Federer-Fleming, Michael-Simon, Almgren.  Late 50's to mid 80's.)

2. Topological dimension theory.  (Work of Brouwer, Lebesgue, 
Szpilrajn.  1900-1940.)

3. Scalar curvature. (Work of Schoen-Yau.  Late 70's.)

4. Hyperbolic geometry and topological complexity.  (Work of Thurston-Milnor.  Late 70's.)

Before turning to the metaphors, I want to discuss why the systolic inequality
is difficult to prove.  

The systolic inequality is reminiscent of the isoperimetric inequality.  Let's
recall the isoperimetric inequality and then compare them.

\newtheorem*{iso}{Isoperimetric inequality}

\begin{iso} Suppose that $U \subset \mathbb{R}^n$ is a bounded open set.  Then
the volume of the boundary $\partial U$ and the volume of $U$ are related by
the formula

$$ Vol_n (U) \le C_n Vol_{n-1}(\partial U)^{\frac{n-1}{n}}. $$

\end{iso}

The isoperimetric inequality is a theorem about all domains $U \subset \mathbb{R}^n$,
and the systolic inequality for the n-torus is a theorem about all the metrics $g$ on
$T^n$.  The set of domains and the set of metrics are of course both infinite.  But in some
practical sense, the set of metrics is larger or at least more confusing.
In my experience, if I make a naive conjecture about all domains $U \subset \mathbb{R}^n$,
with a non-sharp constant $C_n$, the naive conjecture is often right.  If I make a naive conjecture
about all metrics on $T^2$, it is right nearly half the time.  If I make a naive conjecture about
all metrics on $T^3$, it is wrong.  Here's an example.

\newtheorem{nc}{Naive conjecture}

\begin{nc} If $U \subset \mathbb{R}^n$ is a bounded open set, then there is a function $f:
U \rightarrow \mathbb{R}$ so that for every $y \in \mathbb{R}$, the area of the level set
$f^{-1}(y)$ is controlled by the volume of $U$

$$Vol_{n-1}[ f^{-1}(y) ] \le C_n Vol_n (U)^{\frac{n-1}{n}}. $$

\end{nc}

I proved naive conjecture 1 in \cite{Gu1}.

\begin{nc} If $g$ is a metric on $T^2$, then there is a function $f:
T^2 \rightarrow \mathbb{R}$ so that for every $y \in \mathbb{R}$, the length of the level set
$f^{-1}(y)$ is controlled by the area of $g$

$$Length [ f^{-1}(y) ] \le C Area (T^2, g)^{1/2}. $$

\end{nc}

Naive conjecture 2 is also true.  This result is more surprising than the first one.  The problem
was open for a long time.  It was proven by Balacheff and Sabourau in \cite{BS}.

\begin{nc} If $g$ is a metric on $T^3$, then there is a function $f:
T^3 \rightarrow \mathbb{R}$ so that for every $y \in \mathbb{R}$, the area of the level set
$f^{-1}(y)$ is controlled by the volume of $g$

$$Area [ f^{-1}(y) ] \le C Vol (T^3, g)^{2/3}. $$

\end{nc}

Naive conjecture 3 is wrong.  (The counterexamples are based on work of Brooks.   We will discuss
them more in Section 4.)

This anecdote suggests why the systolic inequality is 
much harder to prove in dimension $n \ge 3$.  The basic issue is that the set of 
metrics on $T^3$ is qualitatively larger and stranger than the set of metrics on $T^2$.  In my experience,
the four simple and naive pictures at the beginning of this essay give a fairly decent sample of the possible metrics
on $T^2$.  Let us imagine trying to make a similar sample of metrics on $T^3$.  First of all,
curved three-dimensional surfaces are much harder to picture than curved two-dimensional
surfaces - I don't know how to draw meaningful pictures.  There are metrics on $T^3$ which
are ``analogous" to the two-dimensional pictures above.  But there are also new phenomena
like Brooks's metrics.  These metrics are quite different from any metric on $T^2$, making them
particularly hard to visualize.

Here is another example of a strange high-dimensional metric.

\newtheorem*{gkex}{Gromov-Katz examples}
\begin{gkex} (\cite{K1}) For each $n \ge 2$, and every number $B$, there is a metric on
$S^n \times S^n$ with (2n-dimensional) volume 1, so that every non-contractible n-sphere in $S^n \times
S^n$ has (n-dimensional) volume at least $B$.
\end{gkex}

The Gromov-Katz examples are important in our story, because they show that there is no version of the systolic inequality with n-dimensional
spheres in place of curves.  The first Gromov-Katz examples appear in dimension 4, when $n=2$, but Gromov
and Katz found similar phenomena in dimension 3.

When $n \ge 3$, the zoo of metrics on $T^n$ contains many
wild examples like these.  Because the set of metrics on $T^3$ is so ``big", universal
statements about all the metrics on $T^3$ are rare and significant.

In this essay, I will try to state theorems in the most elementary way that gets across
the main idea.  Therefore, I often don't state the most general version of a theorem.  For example,
Gromov's systolic inequality applies to many manifolds besides tori.  Gromov proved a systolic
inequality for any closed manifold $M$ with $\pi_i(M) =0$ for $i \ge 2$, for real projective spaces,
and for other manifolds.

The four main sections of the paper describe Gromov's four metaphors.  Afterwards
there are two appendices giving other perspectives on the difficulty of the systolic
problem: the lack of good symmetries and the work of Nabutovsky-Weinberger on the 
complexity of the space of metrics.

This essay does not contain an actual proof of the systolic inequality.  For the reader who
would like to learn more, here are some resources.

Gromov's writing on systoles:  The central paper ``Filling Riemannian manifolds" \cite{G1}, Chapter 4 of {\it Metric Structures} \cite{G2},
and the expository essay ``Systoles and isosystolic inqualities" \cite{G4}.

Katz's book {\it Systolic Geometry and Topology} \cite{K2}, and his website on systoles \cite{K3}.

My `Notes on Gromov's systolic inequality' gives in detail Gromov's original proof (14 pages) \cite{Gu4}.

\vskip5pt

{\bf Acknowledgements.} I would like to thank Hugo Parlier for the figure on the first page and Alex
Nabutovsky for helpful comments on a draft of this essay.

\section{The general isoperimetric inequality}

In the late 1950's, Federer and Fleming discovered a version of the isoperimetric inequality for 
n-dimensional surfaces in $\mathbb{R}^N$ for any $n < N$.  This result greatly generalizes
the standard isoperimetric inequality.  Their original result was improved and refined over 25
years until Almgren proved the optimal version in 1986.

\newtheorem*{gii}{General isoperimetric inequality}

Let us write $S^n_R$ to denote a round n-sphere of radius $R$, and $B^n_R$ to denote the 
Euclidean n-ball of radius $R$.

\begin{gii} (Almgren, building on work of Federer-Fleming and Michael-Simon)
Let $M^n \subset \mathbb{R}^N$ be a closed surface with 

$$Vol(M)= Vol(S^n_R). $$

Then there is a surface $Y^{n+1} \subset \mathbb{R}^N$
with $\partial Y = M$ and with 

$$Vol(Y) \le Vol (B^{n+1}_R).$$ 

\end{gii}

Comments. The closed surface $Y$ may not be a manifold.  It will be a manifold
with minor singularities.  The surface $Y$ will always be a chain in the sense of algebraic
topology.  If $M$ is orientable it will be a chain with integer coefficients, and if $M$ is
non-orientable, it will be a chain with mod 2 coefficients.

History.  Federer and Fleming were the first to formulate this inequality \cite{FF}.  In my opinion, 
just formulating
the question was a great contribution to geometry.  The isoperimetric inequality is
the most fundamental and important inequality in geometry.  It's not obvious
how to formulate a version of the isoperimetric inequality for a surface of codimension $\ge 2$.  Such 
a surface does not have a well-defined ``inside".  Instead, Federer and Fleming observed
that a surface $M^n$ of high codimension is the boundary of many surfaces $Y^{n+1}$.  The 
right analogue of the isoperimetric inequality is to claim that one of these many surfaces
has controlled volume.  Federer and Fleming proved the isoperimetric inequality with
a non-sharp constant: $Vol(Y) \le C(n,N) Vol(M)^{\frac{n+1}{n}}$.

In the early 1970's, Michael and Simon improved the constant in this inequality \cite{MS}.  They applied
important ideas from minimal surface theory to the problem, and they were able to prove
that $Vol(Y) \le C(n) Vol(M)^{\frac{n+1}{n}}$.  In other words, their constant $C(n)$ does not
depend on the ambient dimension, and one gets a meaningful inequality for a three manifold
$M$ embedded in some space $\mathbb{R}^N$ of huge or unknown dimension.  In the mid
80's, Almgren proved the sharp constant with a long and difficult proof using geometric measure
theory \cite{A1} .

If $M$ has small volume, then it admits a ``filling" $Y$ whose volume is also small.  The filling $Y$ 
is small in other ways too.  For example, it does not stick out too far away from $M$.  We say
this precisely as follows: we let $N_R(M)$ denote the $R$-neighborhood of $M$.  (In other words, $x \in N_R(M)$
if $dist(x,M) < R$.)  

\newtheorem*{bsfr}{Euclidean filling radius inequality}

\begin{bsfr}(Bombieri-Simon, building on work of Gehring, Federer-Fleming)
Let $M^n \subset \mathbb{R}^N$ be a closed surface with

$$Vol(M) = Vol( S^n_R).$$

Then there is a surface $Y^{n+1}$ with $\partial Y = M$ and with

$$Y \subset N_R(M).$$

\end{bsfr}

Comments. Gromov defined the filling radius of $M \subset \mathbb{R}^N$ as the smallest radius $r$ so
that $M$ bounds some surface $Y \subset N_r(M)$.  According to the filling radius inequality
$Fill Rad(M^n) \le C_n Vol(M)^{1/n}$, with sharp constant $C_n$ coming from the case of a round
sphere.  The filling radius inequality was implicitly proven by Federer and Fleming with a non-sharp constant.  
But Federer and Fleming did not state the inequality.  Gehring formulated his ``link problem" in the 60's -
the link problem is a close cousin of the inequality above.  Gehring proved the inequality with a non-sharp
constant following the method of Federer and Fleming.  Bombieri and Simon proved the sharp inequality
in the early 70's using minimal surface theory \cite{BSi}

Gromov used the Bombieri-Simon inequality to attack the systolic problem for manifolds that embed
nicely into Euclidean space.  If $\Psi: (M^n, g) \rightarrow \mathbb{R}^N$ is a continuous map, we
say that $\Psi$ is an $L$-bilipschitz embedding if, for any two points $p,q \in M$, 

$$\frac{1}{L} | \Psi(p) - \Psi(q) | \le dist_{(M,g)}(p,q) \le L |\Psi(p) - \Psi(q)|. $$

If $(M,g)$ admits an $L$-bilipschitz embedding into Euclidean space, for a reasonably small $L$,
then we can use Euclidean geometry to understand the geometry of $(M,g)$.  In particular, applying
the Bombieri-Simon filling radius estimate, Gromov proved the following inequality:

\newtheorem*{sys1}{Systolic inequality for $(T^n,g)$ nicely embedded in $\mathbb{R}^N$}

\begin{sys1}

If $(T^n,g)$ is a Riemannian n-torus, and there is an $L$-bilipschitz embedding 
from $(T^n, g)$ into $\mathbb{R}^N$, then $(T^n,g)$ contains
a non-contractible curve $\gamma$ with

$$length(\gamma) \le 6 L^2 Vol(T^n,g)^{1/n}.$$

\end{sys1}

Given the Euclidean filling radius inequality, Gromov's proof is about one page long.

At this point, it makes sense to ask whether every $(T^3,g)$ admits an embedding into
some $\mathbb{R}^N$ with bilipschitz constant at most $1000$.  If we had such
embeddings, then we would get the systolic inequality on the 3-torus, and we could
try harder dimensions.  Gromov found that strange metrics $(T^3, g)$ which cannot
be nicely embedded into $\mathbb{R}^N$ no matter how large $N$ is.

\newtheorem*{nex}{Non-embeddable examples}

\begin{nex} (Gromov, 1983) For every number $L$, there is a metric $g$ on $T^3$
so that $(T^3, g)$ does not admit an $L$-bilipschitz embedding into
$\mathbb{R}^N$ for any $N$.
\end{nex}

These examples of Gromov are cousins of the strange examples of Brooks that we mentioned
in the introduction.  In Section 4, we will say a little more about where these examples
come from.

Riemannian manifolds do not admit nice bilipschitz embeddings into Euclidean space.  
But every compact Riemannian manifold $(M^n, g)$ does admit a $1$-bilipschitz 
embedding into $L^\infty$.  In fact, every compact metric space
embeds isometrically into $L^\infty$ as discovered by Kuratowski at the turn
of the century.

\newtheorem*{Ket}{Kuratowski embedding theorem}

\begin{Ket} If $X$ is any compact metric space with distance function $d$,
then there is a map $I$ from $X$ to the Banach space $L^\infty(X)$
so that $d(x,y) = \| I(x) - I(y) \|_{L^\infty}.$
\end{Ket}

To prove the systolic inequality, Gromov extended all the geometric measure
theory described above to the Banach space $L^\infty$.

\newtheorem{met}{Metaphor}
\begin{met} The systolic inequality is like the general isoperimetric inequality
in a Banach space.

\end{met}

This metaphor gives us an approach to the systolic problem.  It doesn't
immediately give us the solution, but it gives an outline of how we may
proceed.  We can take each theorem above and
try to adapt it to the Banach space $L^\infty$.  In this way, we find a lot
of little problems, all related to the systolic inequality, and some of them
easier to approach.

Having said that, the proofs of the general isoperimetric inequality do not
adapt well to the Banach space $L^\infty$.  The Federer-Fleming proof
applies in finite-dimensional Banach spaces, but the constant in their inequality
depends on the ambient dimension, so it doesn't give anything in an infinite-dimensional
space like $L^\infty$.  The proofs of Michael-Simon and Almgren depend
heavily on the Euclidean structure, and they do not adapt to Banach spaces.
The underlying issue seems to be that these proofs exploit the large symmetry group
of Euclidean space.  There are more comments on this issue in Section 5.
Gromov had to rethink the proof of the general isoperimetric inequality - he found
a more robust proof that continues to function in Banach spaces.

The general isoperimetric inequality is a wonderful inequality, and I want to really
encourage people to read about it.  From one point of view, the best proof
is Almgren's proof, because he proves the sharp constant.  But Almgren's proof
is difficult, and it doesn't apply in Banach spaces.  From another point of view,
the best proof is due to Wenger in 2004.  Wenger's proof is only two pages long.
It's very clear, and it needs very few prerequisites.  The proof is as simple
and constructive as the Federer-Fleming argument, the quality of the estimate
is as good as the Michael-Simon argument, and it is even more robust than
Gromov's proof from 1983.  I think anyone working in geometry, analysis,
or topology should find it accessible, and at the same time, it contains a kernel
of wisdom about surface areas which took many years to develop.

\section{Topological Dimension Theory}

In the 1870's, Cantor discovered that $\mathbb{R}^q$ and $\mathbb{R}^n$
have the same cardinality even if $q < n$.  This discovery surprised and
disturbed him.  He and Dedekind formulated the question whether
$\mathbb{R}^q$ and $\mathbb{R}^n$ are homeomorphic for $q < n$.
This question turned out to be quite difficult.  It was settled by Brouwer
in 1909.  Brouwer's theorem is a major achievement of topology.
I'm going to describe some of the history of this result following the
essay ``Emergence of dimension theory" \cite{Jo}.

\newtheorem*{tid}{Topological Invariance of Dimension}

\begin{tid}(Brouwer 1909) If $q < n$, then there is no homeomorphism
from $\mathbb{R}^n$ to $\mathbb{R}^q$.
\end{tid}

Cantor and Dedekind certainly knew that $\mathbb{R}^q$ and $\mathbb{R}^n$
were not {\it linearly} isomorphic.  Linear algebra gives us two stronger
statements:

\newtheorem{linlem}{Linear algebra lemma}

\begin{linlem} If $q < n$, then there is no surjective linear map from $\mathbb{R}^q$ to $\mathbb{R}^n$.
\end{linlem}

\begin{linlem}  If $q < n$, then there is no injective linear map from $\mathbb{R}^n$ to $\mathbb{R}^q$.
\end{linlem}

It seems reasonable to try to prove topological invariance of dimension by generalizing these
lemmas.  A priori, it's not clear which lemma is more promising.
Cantor spent a long time trying to generalize Lemma 1 to continuous maps.
(At one point, Cantor even believed he had succeeded.)  In fact, Lemma 1 does not
generalize to continuous maps.

\newtheorem*{sfc}{Space-filling curve}

\begin{sfc}(Peano, 1890) For any $q < n$, there is a surjective continuous map from $\mathbb{R}^q$ to 
$\mathbb{R}^n$.
\end{sfc}

In his important paper on topological invariance of dimension, Brouwer proved that Lemma 2 does generalize to continuous maps.

\newtheorem*{blt}{Brouwer non-embedding theorem}
\begin{blt} If $n > q$, then
there is no injective continuous map from $\mathbb{R}^n$ to $\mathbb{R}^q$.  
\end{blt}

So it turns out that Lemma 2 is more robust than Lemma 1.  A smaller-dimensional space
may be stretched to cover a higher-dimensional space.  But a higher-dimensional space
may not be squeezed to fit into a lower-dimensional space.  This fact is not
obvious a priori - it is an important piece of acquired wisdom in topology.  In this
section, we're going to talk about the geometric consequences/cousins of this fundamental
discovery of topology.

Shortly after Brouwer, Lebesgue introduced a nice approach to Brouwer's non-embedding
theorem in terms of coverings.  If $U_i$ is an open cover of some set $X$, we say
that the multiplicity of the cover is at most $M$ if each point $x \in X$ is contained in
at most $M$ open sets $U_i$.  We say the diameter of a cover is at most   $\epsilon$ if
each open set $U_i$ has diameter at most $\epsilon$.  For any $\epsilon > 0$, Lebesgue constructed an open
cover of $\mathbb{R}^n$ with multiplicity $\le n+1$ and diameter at most $\epsilon$.  
He then proposed the following lemma.

\newtheorem*{lcl}{Lebesgue covering lemma}
\begin{lcl} If $U_i$ are open sets that cover the unit n-cube, and
each $U_i$ has diameter less than 1, then some point of the n-cube
lies in at least $n+1$ different $U_i$.
\end{lcl}

(Lebesgue proposed his covering lemma in 1909 to give an alternate approach to the
topological invariance of dimension.  In his first paper, he didn't give any proof of the
lemma - perhaps he regarded it as obvious.  Brouwer challenged him to provide a proof,
and a bitter dispute began between the two mathematicians.  Brouwer gave the first
proof of the Lebesgue covering lemma in 1913.)

To see how the Lebesgue covering lemma relates to the non-embedding theorem,
suppose that we have a continuous map $f$ from the unit n-cube to $\mathbb{R}^q$
for some $q < n$.  Lebesgue constructed open covers of $\mathbb{R}^{q}$ with
multiplicity $q+1$ and arbitrarily small diameters.  If $U_i$ cover $\mathbb{R}^q$,
then $f^{-1}(U_i)$ are an open cover of the unit n-cube.  Since $q+1 < n+1$, we
see that some set $f^{-1}(U_i)$ must have diameter at least 1.  On the other hand,
the diameters of the sets $U_i$ are as small as we like.  By taking a limit, we can
find a point $y \in \mathbb{R}^q$ such that the fiber $F^{-1}(y)$ has
diameter at least 1.  So the Lebesgue covering lemma implies the following
large fiber lemma:

\newtheorem*{lfl}{Large fiber lemma}
\begin{lfl} Suppose $q < n$.  If $f$ is a continuous map from the unit n-cube
to $\mathbb{R}^q$, then one of the fibers of $f$ has diameter at least 1.  In other words, there exist
points $p,q$ in the unit n-cube with $|p-q| \ge 1$ and $f(p) = f(q)$.
\end{lfl}

The large fiber lemma is a precise quantitative theorem saying that an n-dimensional
cube cannot be squeezed into a lower-dimensional space.
In particular, the large fiber lemma immediately implies that there is no injective
continuous map from $\mathbb{R}^n$ to $\mathbb{R}^q$.

Gromov thought carefully about the circle of proofs described above, especially
the hypotheses in the Lebesgue covering lemma.  What is it about the unit n-cube
which makes it hard to cover with multiplicity $n$.  Roughly speaking, the key point is that
the unit n-cube is ``fairly big in all n directions", which prevents it from looking like
something lower-dimensional.  Gromov was able to generalize the covering lemma to
spaces that are big in other ways, including spaces with large systole.

\newtheorem*{gfl}{Gromov/Lebesgue covering lemma}
\begin{gfl} (1983) Suppose that $g$ is a Riemannian metric on the n-dimensional
torus $T^n$ with systole at least 10.  In other words, every non-contractible
loop in $(M^n, g)$ has length at least 10.

If $U_i$ is an open cover of $(M^n, g)$ with
diameter at most 1, then some point of $M$ lies in at least n+1 different sets
$U_i$.
\end{gfl}

I haven't looked back at the original papers, but I think that Gromov's proof of 
the covering lemma above extends ideas that originate in Brouwer's original
proof of the covering lemma from 1913.

Topologists following Lebesgue (Menger, Hurewicz...) used the covering
lemma as a basis for defining the dimension of metric spaces.  They said
that the Lebesgue covering dimension of a metric space $X$ is at most $n$
if $X$ admits open covers with multiplicity at most $n+1$ and arbitrarily 
small diameters.  They proved that the Lebesgue covering dimension is
a topological invariant of compact metric spaces.

Different notions of dimension were intensively studied in the first half
of the twentieth century.  The most well-known is the Hausdorff dimension
of a metric space.  The Hausdorff dimension and the Lebesgue covering
dimension may be different.  For example, the Cantor set
has Lebesgue dimension zero and Hausdorff dimension strictly greater than zero.  (The
Hausdorff dimension may be any real number, whereas Lebesgue dimension is always an
integer.)  In 1937, Szpilrajn proved that $Leb Dim(X) \le Haus Dim(X)$ for any compact
metric space $X$.  To do so, he constructed coverings of metric spaces with small
diameters and bounded multiplicities.

\newtheorem*{sdi}{Szpilrajn covering construction}
\begin{sdi} (1937)  If $X$ is a (compact) metric space with n-dimensional Hausdorff measure 0, 
and $\epsilon > 0$ is any number, then there is a covering of $X$ with multiplicity at most
$n$ and diameter at most $\epsilon$.  Hence $X$ has Lebesgue dimension $\le n-1$.
\end{sdi}

Gromov asked whether Szpilrajn's theorem is stable in the following sense: If $X$ has
very small n-dimensional Hausdorff measure, is there a covering of $X$ with multiplicity
at most $n$ and small diameter?

\begin{met} 

The systolic inequality is like a more quantitative version of topological dimension theory - especially Szpilrajn's theorem.

\end{met}

This metaphor gives a second approach to the systolic inequality.  The first half of the approach
is Gromov's systolic version of the Lebesgue covering lemma.  The second half of the approach is a systolic version
of the Szpilrajn theorem which I proved in \cite{Gu2}.

\newtheorem*{rcc}{Covering construction for Riemannian manifolds of small volume} 

\begin{rcc} (Guth 2008) If $(M^n, g)$ is an n-dimensional Riemannian manifold with volume $V$, then there is an open cover of $(M^n, g)$
with multiplicity $n$ and diameter at most $C_n V^{1/n}$.
\end{rcc}

To end this section, we will describe why this covering result is harder than Szpilrajn's, and what kind of new
techniques are needed to prove it.  The key issue is that we need more quantitative estimates.  
The first step in Szpilrajn's proof is to cover our space $X$ with open sets of diameter at most $\epsilon$ and map $X$ to the nerve of the covering.  In Szpilrajn's proof,
we know that the n-dimensional Hausdorff measure of $X$ is zero, and the Hausdorff measure of the image of $X$ is automatically
zero as well.  But in my proof, we know that the n-dimensional volume of $(M,g)$ is some small number $V$.  It does not automatically
follow that the volume of the image of $(M,g)$ is small.  If we choose our cover arbitrarily, then the map to the nerve may stretch
the volume of $(M,g)$ by an uncontrolled factor.  We need to choose an intelligent cover with good estimates on the multiplicity
of the cover, the volumes of the open sets in the cover, the size of the overlaps between neighboring open sets, etc.

Gromov began the job of proving quantitative theorems about open covers in \cite{G3} and \cite{G1}, and my proof builds on 
his ideas.  The main tool is the Vitali covering lemma and variations on it, which give estimates about how
balls can overlap each other.  The Vitali covering lemma first appeared at the beginning of the twentieth century, and it
was used to attack questions in measure theory such as the Lebesgue differentiation theorem.  In the 30's - 50's, it
became an important tool in harmonic analysis, playing a role in the study of the Hardy-Littlewood maximal function,
convolution inequalities, and the Calderon-Zygmund inequalities.  Meanwhile, the covering lemma became an
important tool in geometric measure theory, where it was used to estimate the geometry of surfaces in Euclidean space.
For example, it appears in the Michael-Simon proof of the general isoperimetric inequality, and also more recently in Wenger's
proof of the general isoperimetric inequality.  Gromov began to use the covering lemma to estimate the geometry of
balls in Riemannian manifolds and other metric spaces.

\section{Scalar Curvature}

The scalar curvature is a subtle and important invariant of a Riemannian metric.
It plays an important role in general relativity and also in pure geometry.

The most down-to-earth description of scalar curvature involves volumes of small balls.

\newtheorem*{scb}{Scalar curvature and volumes of balls}
\begin{scb} If $(M^n, g)$ is a Riemannian manifold and $p$ is a point in $M$, then
the volumes of small balls in $M$ obey the following asymptotic:

$$ Vol B(p,r) = \omega_n r^n - c_n Sc(p) r^{n+2} + O(r^{n+3}). \eqno{(*)}$$

\end{scb}

In this equation, $\omega_n$ is the volume of the unit n-ball in Euclidean space,
and $c_n > 0$ is a dimensional constant.  So we see that if $Sc(p) > 0$, then
the volumes of tiny balls $B(p,r)$ are a bit less than Euclidean, and if $Sc(p) <  0$
then the volumes of tiny balls are a bit more than Euclidean.

Understanding the relationship between scalar curvature and the topology of $M$
is a major problem in differential geometry.  Which closed manifolds $M$ admit
metrics with positive scalar curvature?  A guiding problem in the area is the Geroch
conjecture (sadly I am unable to locate the history of this conjecture.  Possibly
I should have attributed it to Kazdan-Warner or to someone else.)

\newtheorem*{ger}{Geroch conjecture}
\begin{ger} The n-torus does not admit a metric of positive scalar curvature.
\end{ger}

The Geroch conjecture was proven by Schoen and Yau in the late 1970's (for $n \le 7$).  Their
proof is one of the main breakthroughs in the study of scalar curvature. 

To get a first sense of the Geroch conjecture, consider the case $n=2$.  In this
case, the scalar curvature is equal to twice the Gauss curvature.  By the Gauss-Bonnet
formula, $\int_{T^2} G darea = 0$.  (Here $G$ denotes the Gauss curvature of
a metric $g$ on $T^2$, and $darea$ denotes the area form of $g$.)  From this formula,
we see that the Geroch conjecture holds for $n=2$.

The Geroch conjecture is much harder for $n \ge 3$.  I would like to try to explain
why.  First of all, the proof we gave for the Geroch conjecture when $n=2$ does not generalize to
higher dimensions.  The Gauss-Bonnet formula does generalize, but the higher-dimensional
version does not involve scalar curvature.

How can we use positive scalar curvature?  Condition $(*)$ about volumes of small
balls sounds comprehensible, but it is quite difficult to apply it.  I guess the key difficulty
is that $(*)$ only applies to the limiting behavior of tiny balls, and it doesn't tell us
anything about balls for any particular radius $r > 0$.  To get a perspective, let's
compare $(*)$ with the Bishop-Gromov inequality for Ricci curvature.

\newtheorem*{bgi}{Bishop-Gromov inequality} 

\begin{bgi} If $(M^n, g)$ is a
Riemannian manifold with Ricci curvature at least 0, then for any
$p \in M$ and any radius $r$,

$$ Vol B(p,r) \le \omega_n r^n. \eqno{(**)}$$
\end{bgi}

The condition $Ric \ge 0$ is much stronger than the condition
$Scal \ge 0$, and the inequality $(**)$ is much stronger than
$(*)$.  Inequality $(*)$ is just a local inequality describing the geometry
of infinitesimal or tiny balls, whereas inequality $(**)$ is a global
inequality, describing the geometry of balls at every scale.  The topology
of a manifold is a global invariant.  It's not so hard to get from a global
geometric estimate like $(**)$ to a theorem about the topology of a manifold,
but there's no way to go immediately from a local estimate like $(*)$ to any
information about the topology or large-scale structure of a manifold.

When $n=2$, the scalar curvature, Ricci curvature, and Gauss curvature
are all equivalent.  In this case, the condition $Scal \ge 0$ implies
$(**)$, and there are plenty of other global geometric
estimates that it implies as well.  But when $n \ge 3$,
the condition $Scal > 0$ definitely does not imply $(**)$.  In fact, it doesn't lead 
to any estimate at all for the volumes of balls of a particular radius, say
$r= 1$.

Bishop, Rauch, Myers, and other geometers proved global geometric inequalities
for manifolds with $Ric \ge 0$ or with $Sec \ge 0$ in the 30's, 40's, 50's.  These
inequalities appeared almost as soon as mathematicians began to look for them.
But a global geometric inequality for metrics with $Scal \ge 0$ was not proven until
the late 1970's, many years after geometers began to look for such an inequality.
It was hard to find partly because you cannot write such an inequality just using
standard geometric quantities like volume, diameter, etc., which appear in the inequalities
of Bishop and Myers.  Instead one has to find new geometric invariants well suited
to the problem at hand.  The first example was the positive mass conjecture.  It
required a lot of wisdom from physics to even formulate the positive mass conjecture.

In the late 70's, Schoen and Yau proved the Geroch conjecture (for dimension $n \le 7$)
as well as the positive mass conjecture.  (See \cite{SY1} and \cite{SY2}.)
The key estimate in their proof is the following observation.

\newtheorem*{syo}{Schoen-Yau observation} 

\begin{syo} If $(M^n, g)$ is a Riemannian manifold
with $Scal > 0$, and $\Sigma^{n-1} \subset M$ is a stable minimal hypersurface,
then $\Sigma$ has - on average - positive scalar curvature also.
\end{syo}

To see how to apply this observation, suppose that $(T^3, g)$ has positive
scalar curvature.  Then a stable minimal hypersurface $\Sigma \subset T^3$ 
is 2-dimensional, and it has (on average) positive scalar curvature.  In two
dimensions, the scalar curvature is much better understood, and it's not so hard to
get topological and geometric information about $\Sigma$.  Now we know topological
and geometric information about every minimal surface $\Sigma$ in $M$, and we can use
this to learn topological and geometric information about $M$ itself.  With this tool, 
Schoen and Yau proved the Geroch conjecture.

Now we can describe Gromov's third metaphor.  As we saw above, the scalar
curvature measures the volumes of tiny (or infinitesimal) balls.  Gromov wondered if
there are similar estimates for the volumes of balls with finite radii.  To make this
metaphor precise, let us define the ``macroscopic scalar curvature" of $(M^n, g)$
at scale $r$ in terms of the volumes of balls with radius $r$.  

Let $p$ be a point in $(M^n, g)$.  
We let $V(p,r)$ be the volume of the ball of radius $r$ around $p$.  Then we let
$\tilde V(p,r)$ be the volume of the ball of radius $r$ around $p$ in the universal
cover of $M$.  (We'll come back in a minute to discuss why it makes sense to use
the universal cover here.)

Now we compare the volumes $\tilde V(p,r)$ with the volumes of balls of radius
$r$ in a constant curvature space.  We let $\tilde V_S(r)$ denote the volume of the
ball of radius $r$ in a simply connected space with constant curvature and scalar
curvature $S$.  Recall that spaces of constant curvature are Euclidean if $S=0$,
round spheres if $S > 0$, or (rescaled) hyperbolic spaces if $S < 0$.  For example,
$\tilde V_0(r) = \omega_n r^n$.  If we fix $r$, then $\tilde V_S(r)$ is a decreasing
function of $S$; as $S \rightarrow + \infty$, $\tilde V_S(r)$ goes to zero, and as
$S \rightarrow - \infty$, $\tilde V_S(r)$ goes to infinity.

If $p \in M$, we
define the macroscopic scalar curvature at scale $r$ at $p$ to be the number $S$ so that
$\tilde V(p, r) = \tilde V_S(r)$.  We denote the macroscopic scalar curvature at scale $r$
at $p$ by $Sc_r(p)$.  In particular, if $\tilde V(p, r)$ is more than
$\omega_n r^n$, then $Sc_r(p) < 0$, and if $\tilde V(p,r) < \omega_n r^n$, then
$Sc_r(p) > 0$.

By formula $(*)$, it's straightforward to check that $\lim_{r \rightarrow 0} Scal_{r}(p) = Scal(p)$.

Let's work out a simple example.  Suppose that $g$ is a flat metric on the n-dimensional
torus $T^n$.  In this case, the universal cover of $(T^n, g)$ is Euclidean space.
Therefore, we have $\tilde V(p,r) = \omega_n r^n$ for each $p \in T^n$ and each
$r > 0$.  Hence $Sc_r(p) = 0$ for every $r$ and $p$.  If we had used volumes of balls
in $(T^n, g)$ instead of in the universal cover, then we would have $Sc_r(p) > 0$
for all $r$ bigger than the diameter of $(T^n, g)$.  By using the universal cover, we
arrange that flat metrics have $Sc_r = 0$ at every scale $r$.

\begin{met} The macroscopic scalar curvature is like the scalar
curvature.
\end{met}

This metaphor leads to some deep, elementary, and wide open
conjectures in Riemannian geometry.

\newtheorem*{ggc}{Generalized Geroch conjecture}

\begin{ggc}(Gromov 1985) Fix $r > 0$.  The n-dimensional torus
does not admit a metric with $Scal_{r} > 0$.
\end{ggc}

The generalized Geroch conjecture is very powerful (if it's true).
Since the scalar curvature is the limit of $Scal_r$ as $r \rightarrow 0$,
the generalized Geroch conjecture implies the original Geroch conjecture.
Taking a fixed value of $r > 0$, the generalized Geroch conjecture implies
the systolic inequality.  Suppose that 
$(T^n, g)$ has systole at least 2.  By the generalized Geroch conjecture,
$Sc_1(p) \le 0$ for some $p \in T^n$.  Therefore, $\tilde V(p, 1) \ge \omega_n$.
Because the systole of $(T^n, g)$ is at least 2, it's not hard to check that
$\tilde V(q,1) = V(q,1)$ for every $q \in T^n$.  Therefore, we see that $V(p,1)$,
the volume of the ball around $p$ of radius 1, is at least $\omega_n$.  Hence
the total volume of $(T^n, g)$ is also at least $\omega_n$.  To summarize,
every metric on $T^n$ with systole $\ge 2$ has volume $\ge \omega_n$.
This is equivalent to the systolic inequality (with a very good constant).

The generalized Geroch conjecture is wide open.  The generalized Geroch
conjecture is considerably stronger than the original Geroch conjecture.
It's also more elementary to state because it only involves the volumes of balls
and not the curvature tensor.  The Geroch conjecture really appeals to me
because it's so strong and so elementary to state, but I don't see any plausible tool
for approaching the problem.  See Appendix 1 for a comment about the difficulty.  

Nevertheless, our third metaphor suggests a different approach to the
systolic inequality, adapting ideas from positive scalar curvature.  In
particular, I was able to adapt the Schoen-Yau estimate for stable
minimal hypersurfaces and prove a weak version of generalized Geroch.

\newtheorem*{wggc}{Non-sharp generalized Geroch}
\begin{wggc} (Guth, 2009) For each $n$, there is a dimensional constant $S(n)$,
so that $T^n$ admits no metric with $Scal_1 > S(n)$.  
\end{wggc}

This result does not give a new proof of the Geroch conjecture.  If we rescale it to
understand $Scal_r$, we see that $\inf_{p \in T^n} Sc_r(p) \le r^{-2} S(n)$.  If
we then take the limit as $r \rightarrow 0$, we get nothing.  Only an absolutely
sharp estimate for $Scal_1$ implies the Geroch conjecture for scalar curvature.
But this result does imply the systolic inequality for the n-dimensional torus.
The constant $S(n)$ works out so that if $(T^n, g)$ has systole at least 2,
then some unit ball in $(T^n, g)$ has volume at least $[8n]^{-n}$.

This technique gives the shortest proof of the systolic inequality for the n-dimensional
torus, but it's not as powerful as other techniques.  For example, recall that a manifold
$M$ is called aspherical if $\pi_i(M) = 0$ for all $i \ge 2$.  There is an old conjecture
that no closed aspherical manifold admits a metric with positive scalar curvature.  The conjecture
is open - the Schoen-Yau technique and other ideas about positive scalar curvature have
not been enough to prove it.  Similarly, this approach to the systolic inequality doesn't
work for all closed aspherical manifolds.  But Gromov's original proof does give the systolic
inequality for all closed aspherical manifolds.  Arguably, Gromov's systolic inequality lends
indirect evidence that aspherical manifolds cannot have positive scalar curvature.

(There are many other techniques in the theory of scalar curvature which may relate to
the systolic inequality.  For example, Gromov and Lawson proved the Geroch conjecture
for all dimensions $n$ in 1979.  Their proof also gives more geometric information
about $(T^n, g)$ than the Schoen-Yau proof.  Their proof is based on Dirac operators,
and there is a key inequality relating the scalar curvature and the spectrum of Dirac operators.
One might ask if there are inequalities relating $Sc_r$ and the spectrum
of the Dirac operator, leading to a Gromov-Lawson approach to the volumes of balls in
Riemannian manifolds.)

\section{Hyperbolic geometry}

Let us return to dimension $n=2$ and consider the systoles of surfaces with
high genus.  There is a systolic inequality for surfaces of high genus, but it
is not as strong as you might expect.  In the 1950's, Besicovitch proved the following
inequality.

\newtheorem*{besi}{Besicovitch systolic inequality}

\begin{besi} If $(\Sigma, g)$ is a closed oriented surface with genus $G \ge 1$, then 

$$ Sys(\Sigma, g) \le \sqrt2 Area(\Sigma, g)^{1/2}. $$

\end{besi}

Let's try to imagine a surface of large genus $G$ with systole around 1.  We
could start with $G$ tori each with systole 1.  Then we could cut out some disks
from the tori, each with circumference around 1, and glue the tori together along
the seams.  If we glue the tori together in a string, then we get a surface of genus
$G$ with systole around 1.  Each of the tori had area $\sim 1$, and so the total
area of our surface is around $G$.  In this way, we get a surface with systole 1
and area around $G$.

It's not at all obvious how to increase the systole of this surface while keeping the area
around $G$.  On the other hand, the systole of this surface is much smaller than
Besicovitch's inequality requires.  The Besicovitch inequality only says that a surface
of area $\sim G$ must have systole at most $\sim \sqrt{G}$.  We will see below that Besicovitch's inequality can be improved a great deal,
but the surface we constructed above can also be improved.  There are surfaces
with genus $G$, area $G$, and systole on the order of $\log G$.  Buser and Sarnak
\cite{BS} gave the first examples of such surfaces: arithmetic hyperbolic surfaces.
These surfaces are among the strangest and most interesting examples in
(Riemannian) geometry.

\newtheorem*{asu}{Arithmetic Hyperbolic Surfaces}

\begin{asu} The isometry group of the hyperbolic plane is $PSL(2, \mathbb{R})$.  If
we take a discrete group $\Gamma \subset PSL(2, \mathbb{R})$, then $\Gamma$
acts on the hyperbolic plane.  Many groups $\Gamma$ act freely, and for such $\Gamma$
the quotient is a hyperbolic surface.

Arithmetic subgroups $\Gamma \subset PSL(2, \mathbb{R})$ lead to particularly
interesting surfaces from the geometric point of view.

For example, define $\Gamma_p \subset PSL(2, \mathbb{Z})$ to be the subgroup
of matrices with modulo $p$ reduction equal to the identity (up to sign).  

$$ \Gamma_p = \left\{  \left( {\begin{array}{cc}  a & b \\   c & d \\ \end{array} } \right) \textrm{ such that } 
 \left( {\begin{array}{cc}  a & b \\   c & d \\ \end{array} } \right)  = \pm  \left( {\begin{array}{cc} 1 & 0 \\ 0 & 1 \\ \end{array}
} \right) \textrm{ modulo } p \right\} .$$

For large prime numbers $p$, $\Gamma_p$ acts freely on the hyperbolic plane.  The resulting
quotient is a non-compact surface with area $\sim p^3$ and genus $\sim p^3$.

This surface is not closed, but there are several tricks for modifying it to get a closed surface.
For example, one can attach small hemispherical caps onto each cusp.  The details
are not important in this essay.

\end{asu}

The arithmetic hyperbolic surfaces play a central role as counterexamples in
Riemannian geometry.  Over the course of this essay, we have mentioned four
strange examples of Riemannian metrics.  Arithmetic hyperbolic surfaces provide
all four strange examples.

1. Large systole.  The surfaces constructed above have genus $G$, area $\sim G$,
and systole $\sim \log G$.  This beats the systole $\sim 1$ for the simple
examples built by gluing together tori.

2. Hard to embed in Euclidean space.  If $(\Sigma, g)$ is a genus $G$ arithmetic hyperbolic surface,
and $\Psi$ is an embedding from $(\Sigma, g)$ to $\mathbb{R}^N$, then the bilipschitz
constant of $\Psi$ is at least $c \log G$.  This estimate does not depend on the dimension
$N$.

Arithmetic hyperbolic surfaces can also be used to construct strange metrics in 
higher dimensions.  For example, let us construct a strange metric on $T^3$.
A surface of any genus may be embedded into $T^3$.  A neighborhood of an
embedded surface $\Sigma$ will be diffeomorphic to $\Sigma \times (-1, 1)$.
On this neighborhood, we can use a product metric, where the metric on $\Sigma$
comes from an arithmetic hyperbolic surface, and the metric on $(-1, 1)$ is the
standard metric with length 2.  Then we extend this metric to the rest
of $T^3$ in such a way that most of the volume is contained in the $\Sigma \times (-1,1)$
region.

3. This metric provides a counterexample to the third naive conjecture in the introduction.  It 
has volume $\sim G$, but if we take any map $F: (T^3, g) \rightarrow \mathbb{R}$, then
one of the level sets will have area at least $\sim G$ also.  Similar constructions
give metrics on $T^3$ that are hard to embed in Euclidean space.

4. Gromov and Katz constructed metrics on $S^n \times S^n$ with large
``n-dimensional systoles", as described in the introduction.  Their original construction
did not use arithmetic hyperbolic surfaces, but later Freedman constructed metrics
with even stronger properties than the Katz-Gromov examples, and Freedman's 
construction is powered by arithmetic hyperbolic surfaces \cite{Fr}.

Arithmetic hyperbolic surfaces are remarkably hard to picture.  When I meet a mathematician
who studies the geometry of surfaces, I often ask them if they have any ideas about visualizing
arithmetic hyperbolic surfaces.  They just laugh.  Part of the problem
is that the systole of an arithmetic hyperbolic surface is only $\sim \log G$.  That means
that to get interesting behavior, we need to look at huge values of $G$.  Naturally,
it is not easy to imagine a surface of genus $10^6$.  Also, many of us try to
visualize Riemannian surfaces as surfaces in three-dimensional Euclidean space.
Arithmetic surfaces embed extremely poorly into Euclidean space, so this strategy
does not work well.

Another possible strategy to get a handle on arithmetic hyperbolic surfaces is to cut
them into simpler pieces.  The most common way to cut a surface into simpler pieces
is called a pants decomposition.  A pair of pants is a surface homeomorphic
to a sphere with three boundary components.  A pants decomposition of a
genus $G$ surface is a set of disjoint simple closed curves on the surface
whose complement
is a union of pairs of pants.  How hard is it to cut an arithmetic surface into
pairs of pants?  For several months, I've been thinking about
how long the curves need to be in a pants decomposition of a genus $G$
arithmetic surface.  Buser constructed a pants decomposition of any genus
$G$ hyperbolic surface using
curves of length $\lesssim G$.  On the other hand, the curves in a pants
decomposition must be larger than the systole $\sim \log G$.  I cannot rule out
a pants decomposition with curves of length $\lesssim \log G$, but I cannot
construct a pants decomposition with curves shorter than $G$.  There
is a tremendous gap between $G$ and $\log G$, and the size of this gap
testifies to my extreme difficulty visualizing arithmetic hyperbolic surfaces.

Arithmetic hyperbolic surfaces are a good example of how algebraic
objects have interesting geometric properties.  See Arnold's essay \cite{Ar} on the
topological efficiency of algebraic objects for more thoughts and perspectives.

Returning to systoles, we have seen that a surface of genus $G$
and area $G$ may have systole around $\log G$.  It turns out that the systole
cannot be bigger than that.

\newtheorem*{hgs}{High genus systolic inequality}
\begin{hgs}(Gromov) Suppose that $\Sigma$
is a closed surface of genus $G \ge 2$ with Riemannian metric $g$.

$$\frac{Sys(\Sigma,g)}{\log G} \le C \left( \frac{Area(\Sigma, g)}{G} \right)^{1/2}. $$

\end{hgs}

This estimate was proven by Gromov in Filling Riemannian manifolds.

So far we've talked about the systoles of some special hyperbolic metrics.  Now I want to go on to a more surprising side of the story: 
using hyperbolic geometry to study the systoles of {\it arbitrary} metrics.  Starting in the 70's
mathematicians used hyperbolic geometry to prove {\it purely topological} theorems about hyperbolic manifolds.  
An important example is the following estimate of Milnor and Thurston.

\newtheorem*{tes}{Triangulation estimate}
\begin{tes}(Milnor-Thurston) Let $(M^n, hyp)$ be a closed hyperbolic
manifold with volume $V$.  Then it requires at least $c_n V$
simplices to triangulate $M$.
\end{tes}

The triangulation estimate has a short, striking proof, which is called the simplex
straightening argument.  See Chapter 5 of \cite{G2} for Gromov's vivid recollection
of learning about the simplex straightening argument.

Philosophically, the Milnor-Thurston result says that a closed hyperbolic manifold 
with large volume is topologically complicated.  Therefore, if you want to build one 
using topologically simple pieces, you will need a lot of pieces.  Here is the analogy 
between the Milnor-Thurston theorem and the systole problem.  Let us suppose that $g$ 
is a metric on the above manifold $M$ with $Sys(M,g) \ge 10$, and consider the unit balls 
in $(M^n, g)$.  It's easy to check that any curve $\gamma$ contained in a 
unit ball $B \subset (M^n, g)$ is contractible.  Now since $(M^n, g)$ has no 
higher homotopy groups, it follows that each unit ball is contractible in $M$.  
Roughly speaking, the condition $Sys(M,g) \ge 10$ forces each unit ball of 
$(M^n, g)$ to be topologically fairly simple.  Since $M$ is topologically complicated, 
it should take a lot of simple pieces to cover $M$, and so one might hope that the volume of $(M^n, g)$ is large.

\begin{met} The systolic inequality for hyperbolic manifolds of large volume is
like the Milnor-Thurston triangulation estimate, but the triangles are replaced by contractible
metric balls.
\end{met}

Following roughly this philosophy, Gromov was able to generalize
the high genus systolic inequality to hyperbolic manifolds of all dimensions.

\newtheorem*{sih}{Systolic inequality for hyperbolic manifolds of large volume} 
\begin{sih} (Gromov 1983) Let $(M^n, hyp)$
be a closed hyperbolic manifold with volume $V > 2$.  Let $g$ be any metric
on $M^n$.  Then $g$ obeys the following systolic inequality:

$$\frac{Sys(M,g)}{\log V} \le C_n \left( \frac{Vol(M, g)}{V} \right)^{1/n}. $$

\end{sih}

A 2-dimensional surface of genus $G \ge 2$ can be given a hyperbolic
metric with volume $2 \pi (2 G - 2) \sim G$.  So as a special case of this theorem, we get
the systolic inequality for 2-dimensional surfaces of high genus stated above.
But Gromov's theorem above applies to surfaces of all dimensions.

The systolic inequality for hyperbolic manifolds is the hardest theorem in systolic
geometry.  The Milnor-Thurston inequality plays a crucial role, but the proof is not
by any means just an adaptation of their proof.  

Unfortunately, there is no expository account of the proof of the hyperbolic systolic
inequality.  Readers may consult Chapter 6.4 of Filling
Riemannian manifolds \cite{G1} or Sabourau's paper \cite{S}.

\section{Appendix 1: Issues of symmetry}

Good mathematical problems often have a lot of symmetry, and
the proofs often find and exploit that symmetry.  Are there any useful
symmetries in the systole problem?  I don't see any obvious symmetries.
In fact, in systolic geometry,
we have to look at some non-symmetric variants of well-known
geometry problems.  The lack of symmetry is one of the main issues
that makes the problems hard.

For example, in Section 1 we saw that the systolic inequality is related
to general isoperimetric inequalities in the Banach space $L^\infty$.
The general isoperimetric inequality in Euclidean space of unbounded
dimension was proven by Michael-Simon and later by Almgren.  One
special feature of Euclidean space is that it's very symmetric, and this
symmetry leads to algebraic formulas that work out nicely.  For example,
the monotonicity formula for minimal surfaces comes from a calculation
that works out nicely in Euclidean space.  The same calculation doesn't
work out nicely in $L^\infty$, and there probably is no monotonicity
formula for minimal surfaces in $L^\infty$.  This monotonicity formula is
the main ingredient in the Michael-Simon proof of the isoperimetric
inequality.  Gromov had to find a different proof, robust enough to work
in non-symmetric spaces.

Here's an example I find even more striking.  Let $PL^\infty$ denote the
projectivization of $L^\infty$.  In other words, $PL^\infty$ is the unit sphere
in $L^\infty$ modulo the action of the antipodal map.  Topologically,
$PL^\infty$ is homotopy equivalent to $\mathbb{RP}^\infty$.  The metric
on $L^\infty$ induces a metric on $PL^\infty$.

We should compare the space $PL^\infty$ and its metric with the standard
metric on real projective space.  Consider the unit sphere in Euclidean
space $\mathbb{R}^{N+1}$.  Take the quotient of the unit sphere by the antipodal map.
The resulting manifold is $\mathbb{RP}^N$ and the resulting metric is
called the Fubini-Study metric.  The Fubini-Study metric is preserved by group
of rotations of $\mathbb{R}^{N+1}$ so it has a large group of symmetries.

The space $PL^\infty$ is much less symmetrical, but it has an important
universal property.  Recall that the Banach space $L^\infty$ has a universal
property: every compact metric space embeds isometrically in $L^\infty$.
Gromov discovered that $PL^\infty$ has an even more striking universal
property.

\newtheorem*{uppl}{Universal property of $PL^\infty$}
\begin{uppl}(Gromov 1983) Suppose that $(\mathbb{RP}^n, g)$ has systole at least 2.  
Then there is a 1-Lipschitz map from $(\mathbb{RP}^n, g)$ into $PL^\infty$, homotopic
to the standard inclusion $\mathbb{RP}^n \subset \mathbb{RP}^\infty$.
\end{uppl}

(Recall that a map is 1-Lipschitz if it decreases all distances.)

Using filling radius techniques, Gromov gave an estimate for the volume of
non-trivial cycles in $PL^\infty$.

\newtheorem*{voc}{Volumes of cycles in $PL^\infty$}
\begin{voc} (Gromov 1983) Any homologically non-trivial n-cycle in $PL^\infty$ must have
volume at least $c_n > 0$.
\end{voc}

Combining the volume estimate and the universal property of $PL^\infty$, we see that
any metric on $\mathbb{RP}^n$ with systole at least 2 has volume at least $c_n > 0$,
which is a systolic inequality for real projective space.  Gromov used his proof
of the systolic inequality to prove this volume estimate.  But if 
one had an independent proof of the volume of cycles estimate in $PL^\infty$, then
we would get a new proof of the systolic inequality for real projective space.

There's a fifth potential metaphor that comes into play here.  In the early 1970's,
Berger and Chern studied the volumes of cycles in $\mathbb{RP}^N$ with the Fubini-Study
metric.  To do so, they used (a variant of) the calibration method invented by
de Rham.

\newtheorem*{calest}{Calibration estimate}
\begin{calest} (Berger-Chern) Let $(\mathbb{RP}^N, g_{FS})$ denote the real projective $N$-space
with the Fubini-Study metric.  If $z^n \subset \mathbb{RP}^N$ is any homologically non-trivial n-cycle, then the volume
of $z$ is at least the volume of a linear copy of $\mathbb{RP}^n \subset \mathbb{RP}^N$.
The lower bound for $Vol(z)$ is one half the volume of the unit n-sphere.  In particular
it does not depend on $N$.
\end{calest}

Gromov's volume estimate for cycles in $PL^\infty$ is analogous to this calibration estimate.  We can
think of this as a fifth metaphor.

\begin{met} The systolic inequality for real projective space is like the Berger-Chern calibration estimate
on $PL^\infty$.
\end{met}

Here is a sketch of the Berger-Chern argument.  Let $P^{N-n}$ be a plane in $\mathbb{RP}^N$ of codimension $n$.
(In other words, $P$ is a linear copy of $\mathbb{RP}^{N-n} \subset \mathbb{RP}^N$.)  
Since $z$ is homologically non-trivial, the topological intersection number of $P$ and $z$ is 1 (mod 2).  Therefore
$z$ and $P$ intersect at least once.  For comparison, let $L$ be
a linear copy of $\mathbb{RP}^n \subset \mathbb{RP}^N$.  The linear space $L$ intersects almost every $N-n$ plane exactly
once.  Hence, for almost every $P$, $z$ intersects $P$ at least as often as $L$ intersects $P$.  Finally, the Crofton
formula tells us that the volume of any n-dimension surface is equal to a fixed constant times the ``average" number of intersections
of the surface with $N-n$ planes $P$.

The Crofton formula is a direct consequence of the symmetry of real projective space (with the Fubini-Study metric).  
To define the ``average" intersection number, we need to define a probability measure on the space of $(N-n)$-planes
in $\mathbb{RP}^N$.  The rotation group acts transitively on the space of planes.  As a compact group, the rotation
group has a natural probability measure  (the Haar measure), which pushes forward to give a unique rotationally-invariant
measure on the space of $(N-n)$-planes.  With this measure, many calculations about ``averages" over the 
$(N-n)$-planes
work out nicely because of the underlying symmetry.  The Crofton formula is a typical example.
Trying to adapt the calibration argument to $PL^\infty$, we meet the lack of symmetry head on, and it's not clear to
me whether the argument can be adapted or not.

There is another side to the symmetry story that I want to mention.  The systole and the other
geometric invariants we have discussed here are extremely robust.  For example, they
are stable under bilipschitz changes of the metric.

\newtheorem*{brs}{Bilipschitz robustness of the systole}
\begin{brs} Let $(M^n, g)$ be a Riemannian manifold.  Let $h$ be another metric on $M$ which 
is L-bilipschitz equivalent to $g$.  (This means that a curve of $g$-length 1 has $h$-length between $L^{-1}$ and $L$.)

Then $L^{-1} Sys(h) \le Sys(g) \le L Sys(h)$.

\end{brs}

This robustness is a kind of approximate symmetry of the systole.  I can make a 3-bilipschitz
change of the metric, and the systole of the new metric will agree with the old systole up to 
a factor of 3.  Let me explain how this is like a symmetry.  Suppose I had an actual symmetry
group for the systole problem.  For each element $\alpha$ of the symmetry group, I could
take any metric $g$ and turn it into a new metric $\alpha g$ with the same systole and the same
volume.  I could try
to use the symmetry group to attack the systole problem as follows: I start with an arbitrary
metric $g$ whose systole I want to estimate.  Then I carefully choose $\alpha$ so that $\alpha g$
is more convenient or more standard than $g$ or has some nice property.  Finally, I estimate
the systole of $\alpha g$.  Now I don't know of any useful group of symmetries for the systole problem.
In other words, I don't know any operations $\alpha$ that I can perform to change the metric $g$
without changing the systole or the volume.  But the bilipschitz robustness of the problem means
that there are lots of changes I can make to $g$ that don't change the systole or the volume very
much.

This robustness is the main symmetry of the systole problem as far as I can see.  Most of
the techniques of systolic geometry aim to exploit it.  They don't exactly use the bilipschitz robustness
stated above, but they use something in a similar spirit.  Because the systole is very robust, 
we can perturb a situation to something simpler and more tractable at the cost of losing a constant factor. 
For this reason, almost every argument in systolic geometry has a non-sharp constant.

Symmetry is very important in mathematics, and this approximate symmetry or robustness is the 
only symmetry I know in the systolic problem.  On the bright side, systolic geometry has
a pretty well-developed system for exploiting this kind of approximate symmetry and proving estimates with non-sharp
constants.  On the dark side, the techniques we have so far offer little for proving sharp
estimates.  For example, we have no idea how to approach the generalized Geroch conjecture.

\section{Appendix 2: Complexity of the space of metrics}

In this section, I want to say a little bit about the work of
Nabutovsky on the complexity of the space
of metrics.  It's remarkable work in metric geometry, and it
gives some insight into the difficulty of proving estimates like
the systolic inequality.

I'm going to discuss a cousin of the systolic inequality:
Berger's isoembolic inequality.  Suppose that $g$ is a metric on $S^n$.  Following
Berger, we define the isoembolic ratio $I(g)$ as follows:

$$I(g) := \frac{Vol (S^n, g)^{1/n}}{Inj Rad (S^n, g)}. $$

The ratio $I(g)$ is scale invariant : it does not change if we rescale
the metric $g$.  Berger proved that the isoembolic ratio is minimized by round
metrics on $S^n$.

\newtheorem*{isoemb}{Isoembolic inequality}
\begin{isoemb} (Berger, 1980, \cite{Be}) Let $g_0$ denote the unit sphere metric on $S^n$.  
Let $g$ denote any metric on $S^n$.

Then $I(g) \ge I(g_0)$.  In other words, if $g$ and $g_0$ have the
same injectivity radius, then $Vol(g) \ge Vol(g_0)$.
\end{isoemb}

We are going to discuss possible ways of proving the
isoembolic inequality.  To set the stage, let's recall the Steiner symmetrization proof
of the classical isoperimetric inequality.
One begins with a domain $U \subset \mathbb{R}^n$, and repeatedly modifies it.
With each modification, the domain becomes more symmetric,
and its isoperimetric ratio improves.  In the limit, the domain converges
to a round ball.  Since the isoperimetric ratio improved with each step,
we can conclude that the isoperimetric ratio of the original domain
was worse than the ratio of the ball, proving the isoperimetric inequality.

This kind of argument appears often in geometry today.  For example,
in Perelman's work on the Ricci flow, one of the minor results is a new proof
of the Gaussian logarithmic Sobolev inequality in Euclidean space.  The
Gaussian logarithmic Sobolev inequality is a variant of the usual Sobolev
inequality, and so it is a cousin of the usual isoperimetric inequality.  For the
purposes of this paper, I think it's most helpful to describe the Gaussian
log Sobolev inequality roughly, leaving out the equations.  The Gaussian
log Sobolev inequality concerns some ratio $S(f)$ where the numerator
is an integral involving $|f|$ and the denominator is an integral involving $| \nabla f|$.
In
Perelman's argument, we begin with an arbitrary non-negative function $f$, and 
we apply a slightly modified heat flow to $f$, giving a family of functions $f_t$.  
Over time, the functions $f_t$ become
more and more symmetric, converging to a standard Gaussian $\gamma$.  Also, the ratio
$S(f_t)$ decreases monotonically.  We then conclude that $S(f)$ is at least
$S(\gamma)$, which is the Gaussian log Sobolev inequality.  (The inequality for non-negative
$f$ implies the inequality for all $f$ by a pretty easy argument.)

It looks tempting to apply this kind of argument to prove the isoembolic
inequality.  The sharp constant in the isoembolic inequality comes from
the round sphere.  Can we begin with a metric $g$ on $S^n$ and
gradually make it more symmetric until it converges to a round
metric while improving the relevant ratio monotonically?  In dimensions
$n \ge 5$, the answer is absolutely not.

\newtheorem*{nwdm}{Nabutovsky huge mountain pass theorem}

\begin{nwdm} (weak version, \cite{N}) Let $n\ge 5$ and let $B$ be any sufficiently large number.
Then there is a metric $g$ on $S^n$ so that $I(g) \le B$, and yet any path $g_t$
from $g$ to a round metric must have $I(g_t) \ge exp ( exp (B))$ for some value of $t$.
\end{nwdm}

The ``mountain pass" refers to the geometry of the graph of the isoembolic ratio $I$.
The metric $g$ given by Nabutovsky's theorem has ``height" $\le B$ - it's moderately
high.  To get from there to the unit sphere metric, one needs to first go up to a very
high height $exp( exp( B) )$, and then come back down - one needs to go over
a huge mountain range separating $g$ from the familiar round metrics.

Let's compare the Gaussian logarithmic Sobolev inequality and the isoembolic
inequality.  The Gaussian logarithmic Sobolev ration $S(f)$ is defined on the
space of functions on $\mathbb{R}^n$.  The isoembolic ratio $I(g)$ is defined
on the space of metrics on $S^n$.  The GL Sobolev inequality says that $S(f)$
attains its minimum at the standard Gaussian.  The isoembolic inequality says
that $I(g)$ attains its minimum at the round metric of any radius.  The Gaussian
is a very symmetrical function, and the round metrics are the most symmetric
metrics on $S^n$.  So far, everything looks similar.  But now, let's move our
attention from the minimizer of the ratio to the graph of the ratio.  Perelman's proof shows
that the
Gaussian logarithmic Sobolev inequality has only one local minimum.  Its
graph looks something like the surface of a parabololoid.  The isoembolic
ratio has infinitely many local near-minima.  These local near-minima can be
extremely deep: to get from a local minimum at height $I$ down to the global
minimum, one may need to first go up to a ridiculous height like $exp( exp( I))$.
The graph of the isoembolic ratio looks something like the craggy rocks on
the bottom of the sea. 

This theorem of Nabutovsky lies near the beginning of a large theory, continued
in joint work of Nabutovsky and Weinberger and still ongoing.  The theory is much stronger and 
more general than what I've presented here.  The function $exp(exp( I(g) ))$ may
be replaced by {\it any computable function of $I(g)$}.  Also, the theorem
applies not just to injectivity radius and volume but to many other
setups.  See \cite{NW} for more information.

However, it is unknown whether the Nabutovsky theorem for injectivity radius
has an analogue for systoles.  We may define a systolic ratio for metrics on $T^n$
as follows:

$$SR(T^n, g) := \frac{Vol(T^n, g)^{1/n}}{Sys(T^n, g))}.$$

The systolic inequality tells us that the infimum of $SR(T^n,g)$ is positive.  We know that
the infimal value is at most $\sim n^{-1/2}$ and is at least $[8n]^{-1}$.  But
besides the minimal value, we know basically nothing about the graph of the function $SR$.  For example,
one has the following open question in the spirit of the work of Nabutovsky
and Weinberger.

\newtheorem*{openques}{Open Question}
\begin{openques} Given a metric on $T^5$ with systole 1 and volume $V$, can it
be deformed to the unit cube metric on $T^5$ while keeping the systole at least 1,
and without increasing its volume too much?
\end{openques}

The Nabutovsky metrics are constructed using logic.  Roughly speaking,
the shape of $(S^n, g)$ encodes an algorithmic problem which is known to be
extremely difficult, and a path $g_t$ from $g$ to the round metric encodes a solution
to the problem.  By logic, one knows that the solution must be very long and complicated.
Nabutovsky is able to reinterpret this geometrically to show that
$I(g_t)$ must be extremely large for some $t$.

\section{Going from two dimensions to three dimensions}

In the first versions of this essay, I didn't mention two dimensions versus
three dimensions.  Then I decided to try to explain why
the systolic inequality is difficult to prove.  The two-dimensional version
was proven in the 40's, and so I tried to say why the three-dimensional
version is harder.  Later, I was trying to give some context about the Geroch
conjecture, so I mentioned how to prove it in two dimensions using Gauss-Bonnet.
Again, the three-dimensional version is much harder.  On the third or fourth draft I noticed a simple thing. 
Almost all of the mathematicians we have been discussing were trying
to generalize a result from two dimensions to higher
dimensions.  

In minimal surface theory, Douglas, Rado, and others solved
the two-dimensional Plateau problem: they proved
the existence of minimal two-dimensional surfaces with prescribed boundary.
Federer and Fleming generalized this result to higher dimensions.
To do so, they invented the general isoperimetric theory described in Section 1.
In the 19th century, Jurgens proved the topological invariance of dimension for 
$\mathbb{R}^1$ and $\mathbb{R}^2$.  Brouwer generalized topological invariance to higher dimensions,
as  discussed
in Section 2.  The two-dimensional version of the Geroch conjecture was proven
by Bonnet in the 19th century.  Schoen and Yau generalized it to higher dimensions, as
discussed in Section 3.  
The possible degrees of 
maps from one (two-dimensional) surface to another were classified by Kneser in the 1930's.  Thurston
and Milnor generalized Kneser's results to three and more dimensions.  To do so, they invented
simplex straightening and proved the triangulation estimate discussed in Section 4.
Loewner proved the systolic inequality in two dimensions in the late 40's.  In the early
80's, Gromov generalized it to dimensions three and higher.  To do so, he invented all of
the metaphors in this essay.  I don't know how Gromov discovered the metaphors in this
paper.  Perhaps he looked for guidance from other mathematicians
who managed to generalize important results from two dimensions to higher dimensions?

Three-dimensional surfaces are far more complicated than two-dimensional surfaces.
Here are some well-known reasons.  First, a curved
three-dimensional surface is much harder to visualize than a curved two-dimensional
surface.  Second, the curvature tensor of a Riemannian three-manifold has six
degrees of freedom at each point, compared to only one degree of freedom for a Riemannian
two-manifold.  Third, there is an explosion in the topological types of closed three-manifolds
as opposed to closed two-manifolds.  Because of the last two points, there are many
strange metrics on three-dimensional manifolds giving counterexamples to naive
conjectures.  Finally, it is very difficult to find ``useful" parametrizations
of Riemannian three-manifolds - whereas the uniformization theorem gives useful parametrizations
for two-manifolds in a wide variety of problems.  In spite of all this complexity, a significant portion of 
two-dimensional geometric theorems remain true in higher dimensions, even when the original
proofs don't generalize.   Can one find some fundamental geometric features that
survive the passage from two dimensions to three dimensions and which tie together
(some of) the subjects discussed in this essay?

To close this essay, let me mention an open problem that marks the edge of my understanding
of metric geometry of three-dimensional surfaces.

\begin{nc} If $g$ is a Riemannian metric on $T^3$, then there is a function $f:
T^3 \rightarrow \mathbb{R}^2$ so that for every $y \in \mathbb{R}^2$, the length of the fiber
$f^{-1}(y)$ is controlled by the volume of $g$

$$Length [ f^{-1}(y) ] \le C Vol (T^3, g)^{1/3}. $$

\end{nc}

This is a naive conjecture that fits into the list of naive conjectures in the introduction.  As a naive
conjecture about three-dimensional metrics, it is probably false, but no one knows yet.  The first place
to look for counterexamples is among the metrics coming from arithmetic hyperbolic surfaces.  I don't know
whether these metrics are counterexamples or not.  To get started, one would need to analyze the lengths
of curves in a pants decomposition of an arithmetic surface, as discussed in Section 4.  There may be
other counterexamples.  Mathematicians have not spent that much time collectively trying to build strange
metrics, and I suspect that many interesting examples are yet to be found.  The work of Nabutovsky
and Weinberger gives some perspective on the difficulty of looking for examples.  

On the other hand, Naive Conjecture 4 may be true.  This conjecture easily implies the systolic
inequality on $T^3$.  It is a much stronger quantitative inequality than what comes from
the methods described in this essay.  Although the question is elementary to state, I don't have any perspective on how to
get started...

\end{document}